\documentclass[12pt]{elsart}
\usepackage{natbib}
\usepackage{amsmath}
\usepackage{amsfonts}
\newtheorem{them}{Theorem}
\newtheorem{lema}{Lemma}
\newtheorem{props}{Proposition}
\newtheorem{rk}{Remark}
\newtheorem{ex}{Example}

\newcommand{\E}{\mathbb{E}}

\begin{document}

\begin{center}
 {\bf \Large Stochastic dynamic programming 
under recursive Epstein-Zin preferences}
\end{center}

\begin{center}  
{\bf  Anna Ja\'skiewicz$^{a}$ and  Andrzej S. Nowak$^{b}$ } 
\end{center}  
\begin{center} \footnotesize{ $^{a}$Faculty of Pure and Applied Mathematics, Wroc{\l}aw University of Science and Technology,  
Wroc{\l}aw, Poland,
{\it email: anna.jaskiewicz@pwr.edu.pl}\\
\noindent $^{b}$Faculty of Mathematics, Computer Science, and Econometrics,  University of Zielona G\'ora,  
Zielona G\'ora, Poland,
  {\it email: a.nowak@im.uz.zgora.pl}} \end{center}

\noindent {\bf Abstract} 
This paper  investigates  discrete-time Markov decision processes  with
 recursive utilities (or payoffs) defined by the classic CES   aggregator and  
 the Kreps-Porteus   certainty equivalent  operator. 
According to the classification introduced by Marinacci and Montrucchio,
some aggregators that we consider are  Thompson and some of them are neither Thompson nor Blackwell. 
We focus on the existence and uniqueness of a solution to the Bellman equation. 
Since the per-period utilities can be unbounded, we work with the weighted supremum norm.
Our paper shows three major points for such models. Firstly,  we prove that  the Bellman equation 
can be obtained by the Banach fixed point theorem for contraction mappings acting on a standard complete metric space. Secondly, 
we  need not  assume any boundary conditions, which are  present when  the Thompson metric or the Du's theorem are used. 
Thirdly, our results give better bounds for the geometric convergence of the value
iteration algorithm 
than those obtained by Du's fixed point theorem. Moreover,  our techniques allow to derive the Bellman equation for 
some values of parameters 
in the CES aggregator and the Kreps-Porteus certainty equivalent that cannot be solved by  Du's theorem
for  increasing and convex or concave operators acting on an ordered Banach space.\\



\section{Introduction}
 
In the classical setting the additive utility  is constructed with the aid of a linear aggregator
of a current utility (or reward) and the expected discounted utility from tomorrow onwards.
More precisely, if $c=(c_n)$ is a consumption sequence and $u$ is the per-period utility, then
the utility $\widehat V_n(c)$  from period $n$ onwards solves the forward-looking equation
$$\widehat V_n(c)=u(c_n)+\beta \E_n(\widehat V_{n+1}(c))\quad\mbox{and}\quad \widehat 
V_n(c)=\E_n\left(\sum_{k\ge n} \beta^{k-n}u(c_k)\right).$$
Here, $\beta\in[0,1)$ is a discount factor and $\E_n$ is the expectation operator with respect to the
time period $n$. Although this framework is mathematically attractable, it does not take into account
the risk associated with transitions from  state to  state.
In order to include the risk aversion into a multi-period decision process, we replace the expected discounted utility
by a certainty equivalent operator induced by an increasing concave utility function defined on the set of
non-negative real numbers. One of the possibilities is to choose a certainty equivalent of the exponential function.
Its negative is sometimes called an entropic risk measure and gains much popularity 
in economics, see \cite{miao,ss}.  \cite{hsar} were among first who applied such a framework to a linear quadratic Gaussian control problem. 
Their  approach  
inspired research   in other areas, e.g.,  in stochastic optimal growth theory
(see   \cite{bjjet}) and in  Markov controlled models  (see \cite{asj}).

Following the seminal work of \cite{kp},
\cite{ez} suggested a recursive utility framework that permits  for a separation 
between risk attitudes and the degree of intertemporal substitution. It 
generalises a deterministic recursive utility introduced by \cite{koop} to a stochastic setting. 
This framework can also allow for early or late
resolution of uncertainty and is tractable to apply because the dynamic programming technique 
can be used.  More precisely, \cite{ez} define 
a recursive utility  under uncertainty by two primitives: a time  aggregator 
and a conditional certainty equivalent. Then, the recursive utility $V_n(c)$ from the period $n$ onwards 
can be written as follows
$$ V_n(c)=\left[(1-\beta)u(c_n)^{1-\rho} + \beta
{\mathcal M}_n(V_{n+1} (c))
^{1-\rho}\right]^{\frac 1{1-\rho}},$$
where $ 0<\rho\not=1.$ Here,
${\mathcal M}_n$ is a  certainty equivalent and  can be induced  for instance by the power function:  
$$ 
{\mathcal M}_n(V_{n+1}(c))=
\E_n(V_{n+1}(c)^{1-\gamma})^{\frac 1{1-\gamma}},\quad 0<\gamma\not=1.
$$
This certainty equivalent is also known under the name  of the Kreps-Porteus operator, see \cite{miao}. 
The value $1/\rho$ represents the 
elasticity of intertemporal substitution between the composite good and the certainty equivalent, 
whereas $\gamma$ governs the level 
of relative risk aversion with respect to atemporal gambles.  Therefore, this time aggregator is known under the name of
CES (Constant Elasticity of Substitution) aggregator.
Recursive preferences   of  this type are called the Epstein-Zin preferences. They form  a core 
of the quantitative asset pricing literature and   find  applications ranging from optimal
taxation to fiscal policy and business cycles, see \cite{by}, \cite{kv} and \cite{sch}.

The growing  popularity of recursive preferences gives rise to 
the expanding domain of their
applications in macroeconomics and finance. The reader is referred in this matter, among others, 
to \cite{backus,miao,ss,sk,x}.
Moreover, the predominant acceptance of recursive preferences  has recently led to a stream of fundamental issues
such as existence, convergence  and uniqueness. For instance, \cite{hs} established a connection 
between the recursive
utility equation and the spectral radius of a corresponding valuation operator. Their approach fostered to
further developments, where \cite{bstach} showed that the spectral radius condition 
is necessary and sufficient for  existence of the Epstein-Zin recursive utility.
Finally,
\cite{chris} making use of this method provided a characterisation of recursive
utility under risk sensitivity, ambiguity and Epstein-Zin preferences. Uniqueness can be obtained
only under  restrictions on the state space and the Markov transition. 
However, none of the aforementioned papers studies dynamic programming with recursive utility.

Despite the substantial progress on recursive utility, 
less is known about its implications for dynamic programming.
Dynamic programming with recursive utility in economics was commenced  by \cite{str} and \cite{os}, who introduced 
a notion of biconvergence. This condition means that utility values can be 
approximated by increasing and
decreasing orbits. Recently, along similar lines, \cite{bich} studied 
deterministic recursive programs under some regularity properties. 
It is worthy to mention that  the work of  \cite{den} was one of the first papers that are related  
to dynamic programming with recursive preferences.

\cite{mm,mmtarski} introduce a contraction-type approach to recursive
utilities. 
More precisely, \cite{mm} divided the Epstein-Zin aggregators
into two classes. In one of them the aggregators are called Blackwell   
and the existence and uniqueness of recursive utility can be proved 
by the Banach contraction mapping principle. 
This approach corresponds to the work of  \cite{blackwell} on discounted dynamic programming.
In the second class, the aggregators are called  Thompson, for which 
the standard methods and metrics cannot be directly applied. The third class of the aggregators 
we deal with are neither Blackwell nor Thompson.
\cite{mm} exploited recursive utilities derived from  the second class and obtained 
significant results by applying the Thompson fixed point theorem
for non-linear operators on ordered normed spaces with a complete and normal positive cone 
(consult with \cite{t}). Another approach  has been recently proposed by 
\cite{brz}, where uniqueness of recursive utility  is presented differently than in \cite{mm,mmtarski}.

Dynamic programming with monotone concave and   $\alpha$-concave  operators  was recently 
adapted by  \cite{bv}  and \cite{balb}, respectively. They demonstrate  new results on the existence 
and uniqueness of a recursive utility as well as some applications to Bellman equations in dynamic optimisation. 
The concave and convex dynamic programming operators have been further developed by \cite{rs1,rs2}.
They obtained valuable results on the existence and uniqueness of solutions of the Bellman equation using
a fixed point theorem of \cite{du}  for monotone convex or concave operators. 
They examined the Epstein-Zin preferences  with unbounded rewards.
However,  this promising approach  based either on the Thompson metric or  Du's theorem 
is frustrated by the fact that uniqueness only obtains in the interior of the domain,
or subject to some appropriate boundary condition. In many relevant applications, this
sort of boundary condition is unnatural. For instance, it imposes restrictions on per-period utilities,
and moreover, it has an influence on the speed of convergence of iterations of dynamic programming operators 
to the value function.
A unified treatment of recursive methods under minimal assumptions on the Koopmans aggregator
has been recently proposed by \cite{bvv}. In particular, they developed an
elementary approach to dynamic programming with recursive utility based on Tarski's
fixed point theorem and establish the existence of a value of the recursive program.

The structure of  non-linear recursive preferences turned out to be harder to study. It was
realised and emphasised in several papers that
the operators involved in the construction of recursive utility are  not contractive in the supremum norm. 
Consequently,
the Epstein-Zin preferences fail to 
satisfy the traditional Blackwell's discounting condition and an appeal
to the Banach contraction mapping principle is of limited scope.  Even establishing
the existence of a recursive utility function has seemed to be an arduous task. 
In this paper, we show that it is not true. 
Our results for Epstein-Zin preferences  are presented within 
an $MDP$ (Markov decision processes) framework.
It worth emphasising that $MDP$ models
play an important  role in operations
research, economics and finance,
see \cite{acemo,br,hl,hll,k,miao,stach}. 

Our work is inspired by the results of  \cite{rs1,rs2}, who  applied  Du's theorem to solve an $MDP$
model with the Epstein-Zin preferences in three cases:
$0<\rho<\gamma<1,$ $ 1<\rho<\gamma$ and    $0<\rho<1<\gamma.$ 
As noted by \cite{by,sch} the case $\rho<\gamma$ is most empirically relevant 
implying that the agent prefers early resolution of uncertainty.
We allow 
the per-period utilities  to be unbounded and use the weighted supremum norms. 
The recursive utility is then a fixed point in a suitably chosen function space.
We  prove the existence and uniqueness of a solution to the Bellman equation in the first two
aforementioned cases  
using the Banach contraction mapping principle for 
the space of  upper semicontinuous functions endowed with the weighted supremum norm, see \cite{w}. 
This method works for some classes with convex aggregators that are according to  \cite{mm} 
of Thompson type. 
Using the Banach fixed point theorem and the supremum norm we get a simple estimation of the speed of  
geometric convergence of the value iteration algorithm to a solution to the Bellman equation. 
As a starting point we can choose a constant function, in one case the zero function. 
We give two examples with bounded rewards showing that 
our value  iterations converge faster than in \cite{rs1,rs2} (that is based on
Du's theorem). Moreover, our regularity 
assumptions on the per-period utility  function and inequalities
involving the weight function 
are weaker than those used in \cite{rs1,rs2}. 
We also show the existence and uniqueness of the solution to the Bellman equation,
with some values of parameters in the Epstein-Zin preferences,
where the dynamic programming operators are neither convex nor concave. 
The fixed point theorem of \cite{du}  in these cases is useless. 
An application of the Banach contraction principle to dynamic programming in $MDPs$ with the Epstein-Zin 
preferences is possible due  to some simple transformations of the Bellman equations.

The paper is organised as follows. In Section  \ref{dp} we provide preliminaries and present our major
theorems. We start with models, where the agent is allowed to use only stationary 
policies.  All proofs from Subsection \ref{F} are collected in Section \ref{proofs}. 
Next in Subsection \ref{NF} we extend our analysis to a non-stationary setting, where history-dependent
policies are permitted. In Section \ref{sec:r} we compare our assumptions imposed on the primitives  
with the ones accepted in \cite{rs1,rs2}. We also provide examples showing that the Banach contraction 
mapping principle gives better bounds for convergence of value iteration algorithm than Du's theorem. 
Finally, Section \ref{sec:c} summarises our main points in this paper.

\section{Dynamic programming model and main results} \label{dp}

\subsection{Preliminaries} \label{preliminary}

Let  $\mathbb{N}$ and $\mathbb{R}$  be the set of all positive integers and all real numbers, respectively. 
Assume that $S$ is a Borel space, i.e., a Borel subset of a complete separable metric space.
We assume $S$  is equipped with its Borel $\sigma$-algebra ${\mathcal B}(S).$
We denote by $\mathcal C$ the Banach space of all bounded continuous functions on $S$ 
endowed with the supremum norm $\|\cdot \|.$
Let
$\omega: S\to [1,\infty)$  be an unbounded continuous function. 
For any $v: S\to\mathbb{R}$ we define  the so-called $\omega$-norm of a function $v$ as follows:
 $$ \|v\|_\omega:=\sup_{s\in S}\frac{|v(s)|}{\omega(s)}.$$
By ${\mathcal C}^\omega$ (${\mathcal U}^\omega$)  we denote the set of all continuous  (upper semicontinuous) functions on $S$ for which 
 $\|\cdot\|_\omega$-norm is finite and by ${\mathcal C}_+^\omega$ (${\mathcal U}_+^\omega$) we denote 
 the set of non-negative functions in ${\mathcal C}^\omega$ (${\mathcal U}^\omega$). 

Assume  that $(X,d)$ is a complete metric space and 
$\delta\in (0,1).$ A mapping 
$F:X\to X$ is called a    $\delta$-contraction,  if 
 $d(Fx_1,Fx_2)\le \delta d(x_1,x_2)$ for all $x_1,x_2\in X.$
By the Banach contraction mapping principle $F$ has a unique fixed point $x_*\in X.$ 
Moreover, for any element $x_0\in X$
\begin{equation}
\label{conv0}
d(F^nx_0,x_*)\le \frac{\delta^n}{1-\delta} d(Fx_0,x_0),
\end{equation}
where $F^n$ is the $n$-th iteration of $F$, see Section 3.12 in \cite{ab}.

Let  $Z$ a separable metric space. By
$\varphi$ we denote a correspondence  from $S$ to $Z.$ For any $D\subset Z$ we
put
$$
\varphi^{-1}(D):=\{s\in S: \varphi(s)\cap D\not= \emptyset\}.
$$
If $\varphi^{-1}(D)$ is a closed (an open) subset in $S$
for each closed (open) subset $D$ of $Z,$ then $\varphi$ is said to be upper
(lower) semicontinuous. If  $\varphi$  is both upper and lower semicotinuous, then it is called continuous.

\subsection{A  model with stationary recursive utilities} \label{ssec:p}

In this paper we study a {\it Markov decision process } 
defined by the following objects.
\begin{itemize}
\item[(i)]  $S$ is a {\it set of states} and is assumed to be a Borel
space.
\item[(ii)]  $A$ is the {\it action space} and
is also assumed to be a Borel space.
\item[(iii)]$\mathbb{K}$ is a non-empty Borel subset of $S\!\times\! A$.
We assume that for each $s\in S$, the non-empty $s$-section $$A(s) =
\{a\in A: (s,a)\in \mathbb{K}\}$$
of $\mathbb{K}$ represents the {\it set of available actions} in state $s$.
\item[(iv)]  $q$ is a Feller transition probability from $\mathbb{K}$
to $S,$ i.e., the function
$$(s,a)\to\int_S v(s')q(ds'|s,a)$$
is continuous for any $v\in{\mathcal C}.$
\item[(v)] $u:\mathbb{K}\to \mathbb{R}$ is a upper semicontinuous non-negative {\it per-period
utility  function}.
\item[(vi)] $\beta\in[0,1)$ is a {\it discount factor.}
\end{itemize}

We first study a {\it stationary model} in which an agent uses the same
decision rule in every period.
Let $\Phi$ be the set of all  Borel measurable
mappings $f$ from $S$ into $A$ such that $f(s)
\in A(s)$ for each $s \in S$. By Theorem 1 in \cite{bp}, $\Phi\not= \emptyset.$  
A {\it stationary policy} is a constant sequence
$\pi = (f,f,\ldots)$ where $f \in \Phi$ and therefore
can be identified with the Borel mapping  $f.$

Assume for a moment that $u$ is bounded.
A {\it state-action aggregator} $H$ maps feasible state-action pairs $(s, a)\in \mathbb{K}$ and bounded Borel measurable 
functions $w$ into real values $H(s, a, w)$ representing lifetime utility, contingent on the current action $a$, the current state $s$ and 
 the use of $w$ to evaluate future states. 
 Traditional additively separable case is implemented by setting:
 $$H(s,a,w)=u(s,a)+\beta \int_S w(s')q(ds'|s,a), $$
where $ (s,a)\in \mathbb{K}.$
Under Epstein-Zin preferences, the  {\it state-action aggregator} $H$  is of the form 
\begin{equation}\label{H}
H(s,a,w) = 
\left[r(s,a)+\beta\left(\int_S w(s')^{  1-\gamma}q(ds'|s,a)\right)^{\frac{1-\rho}{1-\gamma}}\right]^{\frac{1}{ 1-\rho}},
\ \   (s,a)\in \mathbb{K}, 
\end{equation}
assuming that the integral exists
and $0< \rho \not=1, \ 0< \gamma \not=1.$  
Here and in the sequel,
\begin{equation}\label{def:r}
r(s,a)=(1-\beta)u(s,a)^{1-\rho}, \quad (s,a)\in\mathbb{K}.
\end{equation}

\begin{defn}
A function $\widetilde{v}$ is a solution to the {\it Bellman equation}, if 
\begin{equation}\label{bell}
\widetilde{v}(s)=  \max_{a\in A(s)}H(s,a,\widetilde{v})\quad\mbox{for all}\quad s\in S.
\end{equation}
\end{defn}

\begin{defn}
For any policy $f\in \Phi$, we call $v_f$ an $f$-value function, if it satisfies the equation
$$v_f(s)=H(s,f(s),v_f),\quad s\in S$$
with $H$ defined in (\ref{H}).
\end{defn}

The value function $v_f$ is a {\it recursive utility}
in the infinite horizon model in the sense of
\cite{koop}, see also \cite{mm}. It is determined by the aggregator (\ref{H}) and a policy $f\in \Phi.$

\begin{defn}
A stationary policy $f_*\in\Phi$ is {\it optimal}, if
$$ v_{f_*}(s)= \widetilde v(s)\ge v_f(s) \quad\mbox{for all } f\in \Phi\mbox{ and } s\in S,$$
and $\widetilde v$  is a solution to the Bellman equation (\ref{bell}). 
\end{defn}

Such a restriction   to the class of stationary policies and the  recursive utilities in the sense of \cite{koop} is also used by \cite{rs1,rs2}. 

If there exists a policy $f\in \Phi$ that maximises the right-hand side of (\ref{bell}), then, as we show in the sequel,  
$f$ turns out to be an {\it  optimal stationary policy}. This case is well-known in the literature with standard affine
aggregators and expected value as certainty equivalent, see  \cite{bert,hl}. 

\subsection{ Main results on optimal policies in stationary 
models}\label{F}

The  Bellman equation in (\ref{bell})  is precisely of the form  
\begin{equation}
\label{be}
\widetilde v(s) =  \max_{a\in A(s)}  \left[ r(s,a) 
+\beta\left(\int_S \widetilde v(s')^{  1-\gamma}q(ds'|s,a)\right)^{\frac{1-\rho}{1-\gamma}}
\right]^{\frac{1}{ 1-\rho}},  
\end{equation}
for all $s\in S, $ 
provided that  the integral of $\widetilde v(s)^{1-\gamma}$
 exists with respect to $q(\cdot|s,a)$ for every
$(s,a)\in\mathbb{K}.$
Define  $$\theta:= \frac{1-\gamma}{1-\rho}.$$
This number will be used in all four cases of parametrisation in the Epstein-Zin 
preferences considered in this paper.
Let $\hat{0} $ be  the zero function. 
 
In order to find a solution $v$ to the Bellman equation (\ref{be}) when $\gamma, \ \rho \in (0,1)$, we first formulate the following conditions.
\begin{itemize}
\item[(A1)]  For each $s \in  S$, the set $A(s)$
is compact and, moreover, the correspondence  $s\to A(s)$ is upper semicontinuous.
\item[(A2)] There exists a continuous  function $\omega:S\to[1,\infty)$ such that
$$r(s,a)\le M\omega(s), \quad (s,a)\in\mathbb{K}$$
for some constant $M>0.$
\item[(A3)] The function 
$$(s,a) \to \int_S \omega(s')q(ds'|s,a)$$
is  continuous and there is a constant $c\ge 1$ such that
$$    \int_S \omega(s')q(ds'|s,a)\le c\omega(s) $$
 for all $ (s,a)\in{\mathbb{K}}$ and $ c\beta^\theta<1.$
\item[(A4)] The function 
$$(s,a) \to \int_S \omega(s')^\theta q(ds'|s,a)$$
is  continuous on $\mathbb{K}$ and there is a constant $c\ge 1$ such that
$$    \int_S \omega(s')^\theta q(ds'|s,a)\le c\omega(s)^\theta $$
for every $ (s,a)\in {\mathbb{K}} $  and 
$ c^{\frac 1\theta}\beta<1.$
\end{itemize}

\noindent{\bf  Case 1: $0<\rho<\gamma<1.$ }\\
Note that $0<\theta <1$ and for any function $V\in {\mathcal U}_+^\omega$, define the operator $F_1$ as follows:
$$F_1V(s):= \max_{a\in A(s)} H_1 (s,a,V), \ \ s\in S,$$
where
\begin{equation}\label{hatH}
 H_1(s,a,V):=  
 \left[ r(s,a)+\beta\left(\int_S V(s')q(ds'|s,a)\right)^{\frac{1}{\theta}}\right]^\theta,\ \   (s,a)\in\mathbb{K}. 
\end{equation}

We can now state our first main result. 

\begin{them}
\label{thm1} Consider Case 1.
Under assumptions \mbox{(A1)-(A3)}, the following assertions hold.\\
$(a)$ The mapping $F_1:{\mathcal U}_+^\omega \to
{\mathcal U}_+^\omega$ is a $c\beta^\theta$-contraction.
It has a unique fixed point
$w_*.$ Moreover,
if $F^n_1$ is the $n$-th composition of $F_1$ with itself, then 
\begin{equation}
\label{conv1}
\|F_1^n\hat{0} -w_*\|_\omega \le  \frac{c^n \beta^{n\theta}}{1-c \beta^\theta}M.
\end{equation}
$(b)$ The Bellman equation
(\ref{bell}) has a unique solution $ v_* =  w_*^{\frac{1}{1-\gamma}}$.\\
$(c)$ There exists   $f_*\in\Phi$ that satisfies the optimality equation
$$v_*(s)= H(s,f_*(s),v_*)\quad\mbox{for all}\quad s\in S.$$
$(d)$ The stationary policy $f_*$ is  optimal.
\end{them}
 
\noindent{\bf Case 2:}  $0<\gamma<\rho <1.$ \\
Note that $ \theta > 1$ and for any function $V\in {\mathcal U}_+^\omega$, define the operator $F_2$ as follows:
$$F_2V(s):= \max_{a\in A(s)} H_2(s,a,V), \ \ s\in S,$$
where
 $$H_2(s,a,V):=  
 \max_{a\in A(s)} \left[
 r(s,a)+\beta\left(\int_S V(s')^\theta q(ds'|s,a)\right)^{\frac{1}{\theta}}\right] ,\ \   (s,a)\in\mathbb{K}. $$
We assume the integrability of the function $V(\cdot)^\theta$ with respect to $q(\cdot|s,a)$ for $(s,a)\in\mathbb{K}.$  

We can now state our second main result. 
\begin{them}
\label{thm2} Consider Case 2.
Under assumptions \mbox{(A1)-(A2)} and (A4), the following assertions hold.\\
$(a)$ The mapping $F_2:{\mathcal U}_+^\omega \to
{\mathcal U}_+^\omega$ is a $c^{\frac 1\theta}\beta$-contraction.
It has a unique fixed point
$w_*$. Moreover,
if $F^n_2$ is the $n$-th composition of $F_2$ with itself, then 
\begin{equation}
\label{conv2}
\|F_2^n\hat{0} -w_*\|_\omega \le 
   \frac{c^{\frac n\theta} \beta^{n}}{1-c^{\frac 1\theta} \beta}M.
\end{equation}
$(b)$ The Bellman equation
(\ref{bell}) has a unique solution $ v_* =  w_*^{\frac{1}{1-\rho}}$.\\
 $(c)$ There exists   $f_*\in\Phi$ that satisfies the optimality equation
 $$v_*(s)= H(s,f_*(s),v_*)\quad\mbox{for all}\quad s\in S.$$
$(d)$ The stationary policy $f_*$ is  optimal.
\end{them}

If $\gamma, \ \rho>1,$ then we must make some changes
in our assumptions.

\begin{itemize}
\item[(B1)]  For each $s \in  S$, the set $A(s)$
is compact and, moreover, the correspondence $s\to A(s)$ is continuous.
\item[(B2)] There exists a continuous  function $\omega:S\to[1,\infty)$ such that
$$r(s,a)\le M\omega(s), \ \ \mbox{for all} \ \ (s,a)
\in \mathbb{K}$$
and for some constant $M>0.$ The function $r$ is positive and  continuous on $\mathbb{K}.$
\item[(B3)] The function 
$$(s,a) \to \int_S \omega(s')q(ds'|s,a)$$
is  continuous and there is a constant $c\ge 1$ such that
$$    \int_S \omega(s')q(ds'|s,a)\le c\omega(s) $$
 for all $ (s,a)\in{\mathbb{K}}$ and $ c\beta^\theta<1.$
\item[(B4)] The function 
$$(s,a) \to \int_S \omega(s')^\theta q(ds'|s,a)$$
is  continuous on $\mathbb{K}$ and there is a constant $c\ge 1$ such that
$$    \int_S \omega(s')^\theta q(ds'|s,a)\le c\omega(s)^\theta $$
for every $ (s,a)\in {\mathbb{K}} $  and 
$ c^{\frac 1\theta}\beta<1.$
\end{itemize}
The difference between (A3) and (B3) and  between
(A4) and (B4) lies in the properties of the set $\mathbb{K}.$ 

 \noindent{\bf  Case 3: $\gamma> \rho>1.$ }\\
Now we have $\theta >1.$ Therefore, the Bellman equation in (\ref{bell}) has the following equivalent form:
$$w(s)=   \min_{a\in A(s)}
\left[ r(s,a)+\beta\left(\int_S w(s')^\theta q(ds'|s,a)\right)^{\frac{1}{\theta}}\right]$$
for all $s\in S.$
For any function $V\in {\mathcal C}_+^\omega$, define the operator $F_3$ as follows:
$$F_3V(s):= \min_{a\in A(s)} H_3 (s,a,V), \ \ s\in S,$$
where
 $$H_3(s,a,V):=  
 r(s,a)+\beta\left(\int_S V(s')^\theta q(ds'|s,a)\right)^{\frac{1}{\theta}}
$$
 for all $(s,a)\in\mathbb{K}.$ We shall show that $F_3$ maps $ {\mathcal C}_+^\omega $ into itself.
Our next result is as follows.

\begin{them}
\label{thm3} Consider Case 3.
Under assumptions \mbox{(B1)-(B2)} and (B4), the following assertions hold.\\
$(a)$ The mapping $F_3:{\mathcal C}_+^\omega \to
{\mathcal C}_+^\omega$ is a $c^{\frac 1\theta}\beta$-contraction.
It has a unique fixed point
$w_*$  and $ v_* =  w_*^{\frac{1}{1-\rho}}$. Moreover,
if $F^n_3$ is the $n$-th composition of $F_3$ with itself, then 
\begin{equation}
\label{conv3}
\|F_3^n\hat{0} -w_*\|_\omega \le 
\frac{c^{\frac n\theta} \beta^{n}}{1-c \beta}M.
\end{equation}
$(b)$ The Bellman equation
(\ref{bell}) has a unique solution $ v_* =  w_*^{\frac{1}{1-\rho}}$.\\
$(c)$ There exists   $f_*\in\Phi$ that satisfies the optimality equation
$$v_*(s)= H(s,f_*(s),v_*)\quad\mbox{for all}\quad s\in S.$$
$(d)$ The stationary policy $f_*$ is optimal.
\end{them}

 \noindent{\bf  Case 4: $\rho>\gamma >1.$ }\\
Now we have $\theta\in(0,1),$ but $1-\gamma<0$ and $1-\rho<0.$ 
Define the operator $F_4:{\mathcal C}_+^\omega \to  {\mathcal C}_+^\omega $ 
as follows:
$$F_4V(s):= \min_{a\in A(s)}  H_4 (s,a,V), \ \ s\in S,$$
where
$$  H_4(s,a,V):=  
 \left[ r(s,a)+\beta\left(\int_S V(s')q(ds'|s,a)\right)^{\frac{1}{\theta}}\right]^\theta,\ \   (s,a)\in\mathbb{K}. $$
Now we state our last main result in this section. 

\begin{them}
\label{thm4} Consider Case 4.
Under assumptions \mbox{(B1)-(B3)}, the following assertions hold.\\
$(a)$ The mapping $F_4:{\mathcal C}_+^\omega \to
{\mathcal C}_+^\omega$ is a $c\beta^\theta$-contraction.
It has a unique fixed point
$w_*.$  Moreover,
if $F^n_4$ is the $n$-th composition of $F_4$ with itself, then 
\begin{equation}
\label{conv4}
\|F_4^n\hat{0} -w_*\|_\omega \le  \frac{c^n \beta^{n\theta}}{1-c \beta^\theta}M.
\end{equation}
$(b)$ The Bellman equation
(\ref{bell}) has a unique solution  $ v_* =  w_*^{\frac{1}{1-\gamma}}$.\\
$(c)$ There exists   $f_*\in\Phi$ that satisfies the optimality equation
$$v_*(s)= H(s,f_*(s),v_*)\quad\mbox{for all}\quad s\in S.$$
$(d)$ The stationary policy $f_*$ is optimal.
\end{them}
 
The proofs of Theorems \ref{thm1}-\ref{thm4} are given in Section \ref{proofs}.


\begin{rk} { \rm 
$(a)$ If the function  $r$ is bounded, then
in our assumptions we can take $\omega (x)\equiv 1$ and $c=1.$\\
$(b)$ 
The class of functions $r$  and transition probabilities
satisfying our assumptions involving  the weight function $\omega$
is pretty large. For various examples the reader is referred to  \cite{duran} and \cite{jndgaa}. 
Here, we would like to show to construct of an appropriate function $\omega.$
It is quite easy to find a continuous function $\widehat{\omega}:S\to [1,\infty)$ such that
$$r(s,a)\le M \widehat{\omega}(s) \ \ \mbox{for all}\ \  s\in S.$$ 
For any $\beta \in (0,1)$ one can  find $\alpha>1$ such that $\alpha\beta <1.$ Suppose that
the difference 
$$D(s,a):= \int_S\widehat{\omega}(s')q(ds'|s,a)-\alpha\widehat{\omega}(s)>0$$ 
on some subset $\mathbb{K}_0$ of   $\mathbb{K}$ and suppose that $D(s,a)$ is  bounded from above.  
Let $z>0$  and $\omega(s):= z+\widehat{\omega}(s),$
$s\in S.$ Then
$$\int_S\omega(s')q(ds'|s,a) = \alpha \omega(s)
+ (1-\alpha)z +D(s,a).$$ 
If we take $z$   such that $(1-\alpha)z +D(s,a)<0$ for 
all $(s,a)\in \mathbb{K}_0,$ 
  then
$$\int_S\omega)s')q(ds'|s,a) \le  \alpha \omega(s)\ \
\mbox{for all}\ \ (s,a)\in\mathbb{K}. 
$$
In this way, we can establish an appropriate  weight function 
for several  cases arising in economic dynamics.
 }
\end{rk}

\begin{rk}\label{rem2} {\rm 
From (B2) it follows that  
$$r_0(s):=\min_{a\in A(s)}r(s,a)>0\ \ \mbox{ for all}\ \ s\in S.$$ 
Hence, the operators $F_3$ and $F_4$ are  considered 
on ${\mathcal C}_+^\omega.$ In particular, $F_3\hat{0}$ and 
$F_4\hat{0}$ are well-defined. Therefore, the fixed point  of $F_i$ (for 
$i=3,4$) can be obtained by  iterating the operator on the
zero function, which makes the  approximation procedure simpler.}
\end{rk}

\begin{rk} {\rm Cases 1 and 3 were considered in \cite{rs1,rs2}. Their proof of the existence of a unique solution 
to the Bellman equation is based on the fixed point theorem of \cite{du} for monotone convex or concave operators 
on ordered Banach spaces. They also show the geometric convergence of iterations of the corresponding operators.
However, this speed of convergence  is slower than in  Theorems \ref{thm1} and \ref{thm3}. 
We discuss this issue  in Subsection \ref{gc}. 
The results of \cite{rs1,rs2} also apply to the case when $0<\gamma <1<\rho$ which is not covered by our approach. 
On the other hand, Cases 2 and 4 
cannot be solved by using the theorem of \cite{du}, since the corresponding operators  are in these aforementioned cases
neither concave nor convex. For more comments on this issue
see Subsection \ref{duthm}. }
\end{rk}

\begin{rk} {\rm Note that our conditions in Theorem \ref{thm3} are weaker than those in \cite{rs1}, who assumed that
the  reward function is bounded from below by some function greater than 1.
Moreover, they also impose that the integral of the weight function is bounded from below by this weight function
multiplied by some constant smaller than $1/\beta^\theta.$
The reader is referred to Assumption 4.1   and inequalities  (42) and (44) in \cite{rs1}.}
\end{rk}
 
\begin{rk} 
{\rm In the most empirically relevant parametric specifications  $\rho<\gamma$,
implying that the agent prefers early resolution of uncertainty (Cases 1 and 3). 
For a discussion of this topic the reader is referred to  \cite{by} or \cite{sch}.  
Cases 2 and 4 where 
the coefficient of relative risk aversion $\gamma >1$ and 
the intertemporal elasticity of substitution $\frac{1}{\rho} <1$ 
are also relevant, see 
\cite{by,basu,fw}.}
\end{rk}

\subsection{An extension to models with non-stationary policies} \label{NF}
 
A definition  of the decision process involving  the Epstein-Zin aggregators in a framework of 
general class of history dependent policies
is not obvious. For some specifications of the parameters 
$\gamma$ and $\rho$, one can define utilities  in finite step models and then consider 
the limit as the time horizon tends to infinity. 
Then using   iterations of some appropriately defined  dynamic programming operators,   
one can  prove that an optimal stationary policy derived from
the Bellman equation is also optimal in the class  of all policies. 

A similar construction in  risk-sensitive Markov decision processes involving the additive (affine)  aggregator and 
entropic risk measure as a certainty equivalent
is described in  \cite{asj,bjjet}. Additional comments on this issue can be found in \cite{ss}.

A   {\it policy} is a   sequence of Borel measurable mappings
$\pi = (\pi_{n})_{n\in\mathbb{N}}$ from the history space to the action set. More precisely,
each $\pi_{n}(h_{n})\in A(s_n)$, $n\in\mathbb{N},$ where
$h_{n} = (s_{1},a_{1},\ldots,s_{n-1},a_{n-1},s_{n})$ is the history of the process
up to the $n$-th state. For $n=1,$ $h_1=s_1.$ 
The {\it class of all policies} is denoted by $\Pi$.
A {\it Markov policy} is a   sequence  
$\pi = (\pi_{n})_{n\in\mathbb{N}}$ where each $\pi_n$ depends only on the $n$-th state $s_n\in S.$ The   
{\it class of all Markov policies} is denoted by $\Pi_M$.

To avoid  complex  notation, we shall confine ourselves to consideration of  {\it Markov policies}. 
The proof that Markov policies are enough in dynamic programming is  standard  and will be omitted in this paper.
With respect to this issue consult with \cite{bert,hl,ss}.

Let us fix a Markov policy $\pi = (\pi_{n})_{n\in\mathbb{N}}$.  For a non-negative function $\psi: S\to\mathbb{R}$
define the operator 
\begin{eqnarray}\label{Top}
{\mathcal T}_{\pi_n}\psi(s_n) := 
\left[r(s_n,\pi_n(s_n))+\beta\left(\int_S \psi(s_{n+1})^{ 1-\gamma} 
q(ds_{n+1}|s_n,\pi_n(s_n))\right)^{\frac{1-\rho}{1-\gamma}}\right]^{\frac{1}{ 1-\rho}}
\end{eqnarray}
 for all $s_n\in S.$ We assume again the integrability of the function $\psi(\cdot)^{1-\gamma} $
 with respect to $q(\cdot|s,a)$ for every $(s,a)\in\mathbb{K}.$ If $1-\gamma<0,$ then $\psi$ must be positive. 
 Recall that $r$ is defined in (\ref{def:r}).
The operator in (\ref{Top}) is monotone
for some values of parameters $\gamma$ and $\rho$. Here, we shall consider cases  discussed in the previous section.

\noindent{\bf Case (A)}: $\gamma,\ \rho \in (0,1).$ Moreover,
(A1)-(A2) are satisfied and  if $\theta\in(0,1)$, then  (A3) holds and if $\theta>1,$ then (A4) holds.

The recursive utility in the finite time horizon is defined by the following composition of the operators
\begin{equation}
\label{utn}
 U_{n}(\pi)(s_1):= 
{\mathcal T}_{\pi_1}{\mathcal T}_{\pi_2}\cdots {\mathcal T}_{\pi_n}\hat{0}(s_1).
\end{equation}
Here, the  process starts at  state $s_1 \in S$  and terminates after $n$ steps.
with the terminal utility which equals zero. 
Since $\gamma,\ \rho \in (0,1),$ it follows that 
$U_{m}(\pi)(s_1) \le U_{m+1}(\pi)(s_1)$ for all $m\ge 1.$ 
The   {\it  recursive utility} of an agent in the infinite time horizon  
 is defined as
$$U(\pi)(s_1):= \lim_{m\to\infty}U_m(\pi)(s_1), \ \  s_1\in S,\  \pi\in \Pi_M.$$ 

\begin{defn}\label{defut2}
A policy $\pi_* $ is optimal, if
$$U(\pi_*)(s_1)  \ge U(\pi)(s_1)\quad\mbox{for all}\quad \pi\in \Pi,\ s_1\in S.$$ 
\end{defn}

The composition ${\mathcal T}_{\pi_1}{\mathcal T}_{\pi_2}\cdots {\mathcal T}_{\pi_n}\psi(s_1)  $
can be defined for any $\psi\ge 0$ and we have
$$
{\mathcal T}_{\pi_1}{\mathcal T}_{\pi_2}\cdots {\mathcal T}_{\pi_n}\psi(s_1)  \ge {\mathcal T}_{\pi_1}{\mathcal T}_{\pi_2}\cdots {\mathcal T}_{\pi_n}\hat{0}(s_1).
$$
Moreover, if $\pi$ is a stationary policy $(f,f,...,)$ identified
with $f\in \Phi,$ then from the above definitions and Lebesgue monotone
convergence theorem, it follows that
\begin{equation}
\label{fpp1}
U(f)(s_1)= H(s_1,f(s_1), U(f))\ \ \mbox{for all }\ \ s_1\in S,
\end{equation}
where $H$ is given in (\ref{H}). Let us define 
\begin{equation}
\label{tmax}
{\mathcal T}\psi(s_n):= \max_{a_n\in A(s_n)}
H(s_n,a_n,\psi),\ \ \psi\ge 0.
\end{equation}
 Note that
the sequence ${\mathcal T}^m\hat{0}$ of $m$ compositions of the operator ${\mathcal T}$ with itself
is non-decreasing and for any Markov policy $\pi$, we have
$$U(\pi)(s_1) = \lim_{m\to\infty}
{\mathcal T}_{\pi_1}{\mathcal T}_{\pi_2}\cdots {\mathcal T}_{\pi_m}\hat{0}(s_1) \le
 \lim_{m\to\infty}{\mathcal T}^m\hat{0}(s_1)\ \ \mbox{for all}\ \ s_1\in S.
$$

\begin{them}\label{thm5} In Case (A), there exists a stationary
optimal policy in the class of all policies. 
\end{them}

\noindent{\bf Proof  }
By Theorem \ref{thm1} or  \ref{thm2} (depending on whether (A3) or (A4) is assumed),
the Bellman equation has a unique solution $v_*$ (in appropriate function space) and there exists 
$f_*\in \Phi$ such that
$$ v_*(s_1)= {\mathcal T}v_*(s_1) = H(s_1,f_*(s_1),v_*)\quad\mbox{for all}\quad s_1\in S.$$
By (\ref{fpp1}) we also have 
$U(f_*)(s_1)= H(s_1,f_*(s_1), U(f_*))$ for   $s_1\in S.$
This fact and the proof of Theorem  \ref{thm1} or  \ref{thm2}   imply that 
$v_*(s)=U(f_*)(s)$ for all $s\in S.$  Let   $\pi = (\pi_{n})_{n\in\mathbb{N}}$ be any Markov policy. 
We have $v_*(s)={\mathcal T}v_*(s)\ge {\mathcal T}_{\pi_n}v_*(s)$ for all $n \in \mathbb{N}$ and $s\in S.$ Iterating these
inequalities we obtain for any $m\in\mathbb{N}$ that
\begin{eqnarray*}
v_*(s)&\ge& {\mathcal T}_{\pi_1}v_*(s_1) \ge 
{\mathcal T}_{\pi_1}{\mathcal T}_{\pi_2}v_*(s_1) \ge \cdots \ge 
{\mathcal T}_{\pi_1}{\mathcal T}_{\pi_2}\cdots {\mathcal T}_{\pi_m}v_*(s_1) \\
&\ge& {\mathcal T}_{\pi_1}{\mathcal T}_{\pi_2}\cdots {\mathcal T}_{\pi_m}\hat{0}(s_1)=
U_m(\pi)(s_1).  \end{eqnarray*} 
Hence, it follows that
$$U(f_*)(s_1)=v_*(s_1)\ge \lim_{m\to\infty}U_m(\pi)(s_1)= U(\pi)(s_1)$$  
for every $s_1\in S.$  $\Box$
	
\noindent{\bf Case (B)}: $\gamma>1$ and $ \rho >1.$
Moreover, (B1)-(B2)  are satisfied and  if $\theta\in(0,1),$ then (B3) holds and if  $\theta>1,$ then  (B4) holds.

For any $f\in\Phi$ we define  ${\mathcal T}_f\hat{0}(s):=
r(s,f(s))^{\frac{1}{1-\rho}}$ for $s\in S.$  For $\psi>0$  we put
${\mathcal T}_f\psi(s):= H(s,f(s),\psi).$ Then 
for any Markov policy  $\pi = (\pi_{n})_{n\in\mathbb{N}}\in \Pi_M$,
the recursive utility  in the finite time horizon $U_n(\pi)(s_1)$ can be defined as in
(\ref{utn}), i.e.,
$$ U_n(\pi)(s_1) = 
{\mathcal T}_{\pi_1}{\mathcal T}_{\pi_2}\cdots {\mathcal T}_{\pi_n}\hat{0}(s_1),\ \  n\in\mathbb{N}.$$
 It should be noted that if $1-\gamma<0$ and $1-\rho<0,$ then we  have $U_{m+1}(\pi)\le U_m(\pi).$
Therefore, we can define
$$ U(\pi)(s_1):= \lim_{m\to\infty}U_m(\pi)(s_1).$$
The sequence ${\mathcal T}^m\hat{0}$ of $m$ iterations  of the  operator ${\mathcal T}$ defined in (\ref{tmax}) with itself  is non-increasing.
Under our assumptions
${\mathcal T}^n\hat{0}(\cdot)$ is continuous on $S$ for every  $n\in\mathbb{N}.$ Let 
$$w_0(s):= \lim_{m\to\infty}{\mathcal T}^m\hat{0}(s) \quad \mbox{for  }s\in S.$$ 
The function $w_0$ is upper semicontinuous and using Proposition 10.1 in \cite{schal}, we infer that for all $s\in S$ it holds
\begin{eqnarray}
\label{w0}\nonumber
w_0(s)&=& \lim_{m\to\infty}{\mathcal T}^m\hat{0}(s) = \lim_{m\to\infty} \max_{a\in A(s)} H(s,a,{\mathcal T}^{m-1}\hat{0})=
 \max_{a\in A(s)}\lim_{m\to\infty}  H(s,a,{\mathcal T}^{m-1}\hat{0}) \\&=&\max_{a\in A(s)} H(s,a,w_0).
\end{eqnarray} 

\begin{them}
\label{thm6}
In Case (B), there exists a stationary optimal policy in the class of all policies. 
\end{them}

\noindent{\bf Proof  } 
By Theorem   \ref{thm3} or \ref{thm4} (depending on whether (B3) or (B4) is assumed),
the Bellman equation has a unique solution $v_*$ and there exists  $f_*\in \Phi$ such that
$$ v_*(s_1)= {\mathcal T}v_*(s_1) = H(s_1,f_*(s_1),v_*)\quad\mbox{for all}\quad s_1\in S.$$
Our arguments used in the proofs of Theorems 
   \ref{thm3} and \ref{thm4}  and (\ref{w0}) imply that $w_0(s)=v_*(s)=
U(f_*)(s)$ for all $s\in S.$ Fix any state $s_1\in S$, a number $\epsilon>0$ and 
a Markov policy   $\pi = (\pi_{n})_{n\in\mathbb{N}}.$  
There exists $N>0$ such that
$$w_0(s_1) =v_*(s_1) \ge {\mathcal T}^m\hat{0}(s_1) -\epsilon 
\ge U_m(\pi)(s_1)-\epsilon,\quad\mbox{for all}\quad m >N.$$
Hence
$$
U(f_*)(s_1)=v_*(s_1) \ge \lim_{m\to\infty}
( U_m(\pi)(s_1) -\epsilon)= 
U(\pi)(s_1) -\epsilon.$$
Since $\epsilon >0$ can be arbitrarily small,
it follows that $U(f_*)(s_1) \ge U(\pi)(s_1).$  $\Box$

\section{Proofs  of Theorems \ref{thm1}-\ref{thm4}}\label{proofs}

The proof of Theorem \ref{thm1} is proceeded by an auxiliary lemma.

\begin{lema}
\label{l1}
Under assumptions \mbox{(A1)-(A3)}, the  assertions $(a)$-$(b)$ of Theorem \ref{thm1} hold. 
\end{lema}

\noindent{\bf Proof}\ \   First, we must show that 
 $F_1$ maps ${\mathcal U}_+^\omega$ into itself.   Let $V\in {\mathcal U}_+^\omega.$
Note that from  (A2) and (A3), we get
\begin{eqnarray*} 
0&\le & F_1V(s) \le  \left[M\omega(s) +\beta \left(\int_S \|V\|_\omega\omega(s')q(ds'|s,a)\right)^{\frac{1}{\theta}}\right]^\theta
\\ &\le &
 \left[M\omega(s)^{\frac 1\theta} +\beta \left( c\|V\|_\omega \omega(s)\right)^{\frac{1}{\theta}}\right]^\theta \le 
 \omega (s)\left[ M+\beta c^{\frac 1\theta}\|V\|_\omega^{\frac 1\theta}\right]^\theta,\quad s\in S.
\end{eqnarray*}
Additionally, by Lemma 8.5.5 in \cite{hll}, it follows that the function 
$$(s,a) \to \int_S V(s')q(ds'|s,a)$$
is upper  semicontinuous. 
Therefore, $(s,a)\to H_1(s,a,V)$ is upper semicontinuous. Consequently, by (A1) and Lemma 3.8.3(c) in \cite{hll}
the function  $F_1V$ is upper semicontinuous on $S.$ 
Next we prove that $F_1$ is a $c\beta^\theta$-contraction. 
Let $V_1, V_2\in {\mathcal U}^\omega_+$. Consider  the norm 
$$\|\pmb{z}\|_{\frac{1}{\theta}}=
\left(|z_1|^{\frac{1}{\theta}}+ |z_2|^{\frac{1}{\theta}}\right)^\theta$$
for any vector $\pmb{z}=(z_1,z_2)\in{\mathbb{R}}^2.$
 Observe that $\|\cdot\|_{\frac1\theta}$ is indeed the norm, since $\frac 1\theta>1.$
Now take
$$\pmb{x}:=\left(r^\theta,\beta^\theta t_1\right), \quad  \pmb{y}:=\left(r^\theta,\beta^\theta t_2\right),$$
where we use the following abbreviations
$$r^\theta:=r(s,a)^\theta, \quad t_i:=\int_S V_i(s')q(ds'|s,a)$$
for $(s,a)\in\mathbb{K}.$
We notice that 
$$\|\pmb{x}\|_{\frac 1\theta}=H_1(s,a,V_1) \ \mbox{ and }\  \|\pmb{y}\|_{\frac1\theta}=  H_1(s,a,V_2)$$
with $H_1$ in (\ref{hatH}). Although for convenience of notation we do not indicate dependence of $\pmb{x}$ and 
$\pmb{y}$ on $(s,a)$ we keep it in mind.  Using these facts, we obtain
\begin{eqnarray*}
|F_1V_1(s)-F_1V_2(s)|  &\le& \max_{a\in A(s)} | H_1(s,a,V_1)- H_1(s,a,V_2)| =
\max_{a\in A(s)}\left|\|\pmb{x}\|_{\frac 1\theta}-\|\pmb{y}\|_{\frac 1\theta}\right| \\
&\le&
\max_{a\in A(s)}   \| \pmb{x}-\pmb{y} \|_{\frac 1\theta} 
=\max_{a\in A(s)} \|(0,\beta^\theta(t_1-t_2))\|_{\frac 1\theta} \\
&=&\max_{a\in A(s)} \left( 0+\left|\beta^\theta (t_1-t_2) \right|^{\frac 1\theta} \right)^\theta
= \max_{a\in A(s)}\beta^\theta |t_1-t_2|
 \\  &\le &
\max_{a\in A(s)} \beta^\theta  \int_S \frac{|V_1(s')-V_2(s')|}{\omega(s')} \omega(s')q(ds'|s,a)\\
&\le& \beta^\theta c \|V_1-V_2\|_\omega \omega(s), \quad s\in S.
\end{eqnarray*}
Now dividing all sides of these inequalities by $\omega(s)$ and taking the supremum over all $s\in S$  we get
$$\|F_1V_1-F_1V_2\|_\omega \le c \beta^\theta  \|V_1-V_2\|_\omega.$$
Hence, $F_1$ is a $c\beta^\theta$-contraction. 
Since ${\mathcal U}_+^\omega$ is a closed subset of the Banach space ${\mathcal C}^\omega$ we conclude
from the Banach contraction principle  that there a unique $w_*\in{\mathcal U}^\omega_+$ such that
$$w_*(s)=F_1w_*(s)\quad\mbox{for } s\in S.$$
Putting $v_*:=w_*^{\frac 1{1-\gamma}}$ we get that $v_*$ is non-negative, upper semicontinuous and 
\begin{equation}\label{be1}
 v_*(s) =\max_{a\in A(s)}  
  \left[ r(s,a)+\beta\left(\int_S  v_*(s')^{  1-\gamma}q(ds'|s,a)\right)^{\frac{1-\rho}{1-\gamma}}
\right]^{\frac{1}{ 1-\rho}} 
\end{equation}
for all $s\in S.$ To finish this proof note that (\ref{conv1}) now follows from (\ref{conv0}) with $F=F_1.$  $\Box$

\noindent{\bf Proof of Theorem \ref{thm1}}\ \ In view of Lemma \ref{l1} it remains to prove assertions $(c)$ and $(d).$
The existence of a maximiser $f_*\in \Phi$ of the right-hand side of (\ref{be1}) follows from our assumptions (A1)-(A3)
and Corollary 1 in  \cite{bp}. Let ${\mathcal B}^\omega$  be the set of all Borel measurable functions on $S$ for which 
 $\|\cdot\|_\omega$-norm is finite and  ${\mathcal B}_+^\omega$ be
 the set of non-negative functions in ${\mathcal B}^\omega.$ For any $f\in \Phi$ and 
 $w\in{\mathcal B}_+^\omega$, we define the operator $F_{1f}$ as follows
 $$F_{1f}w(s):=H_1 (s,f(s),w),\quad s\in S.$$
Using the measurability results of \cite{bp}, similarly as above, we can show that $F_{1f}w\in {\mathcal B}_+^\omega.$ 
Moreover, $F_{1f}$ is a $c\beta^\theta$-contraction.
Therefore, the operator $F_{1f}$ has s a unique fixed point $w_f\in{\mathcal B}_+^\omega$ and for any function $w_0\in  {\mathcal B}_+^\omega$
$$\lim_{n\to\infty} \|F^n_{1f}w_0-w_f\|_\omega=0.$$
Observe also that $F_{1f}$ is monotone, i.e. $F_{1f}w_1\ge F_{1f}w_2$ whenever $w_1\ge w_2$ and $w_1,w_2\in{\mathcal B}_+^\omega.$  
For any $f\in \Phi$, we have $F_{1f}w_*\le w_*.$ By iteration of this inequality we get
$F^n_{1f}w_*\le w_*$ and consequently $w_f\le w_*,$ since $w_*\in{\mathcal U}_+^\omega\subset {\mathcal B}_+^\omega.$
On the other hand,  due to Lemma \ref{l1} we have  $F_{1f_*}w_*=w_*=F_1w_*.$ 
 Again iterating $F_{1f_*}$ we obtain $\lim_{n\to\infty} F^n_{1f_*} w_*=w_{f_*}=w_*.$
Summing up, 
$$
w_f\le w_*=w_{f_*}\quad\mbox{for all }\ f\in \Phi,
$$
 which implies via transformation $v:=w^{\frac 1{1-\gamma}}$ with $1-\gamma>0$ that
$$
v_f\le v_*=v_{f_*}\quad\mbox{for all }\ f\in \Phi.
$$
This means that  $f_*$ is  an optimal stationary policy.
 $\Box$


The proof of Theorem \ref{thm2} is proceeded by an auxiliary fact.

\begin{lema}
\label{l2}
Under assumptions (A1)-(A2) and (A4), the  assertions $(a)$-$(b)$ of Theorem \ref{thm2} hold. 
\end{lema}

\noindent{\bf Proof}\ \   First, we must show that 
 $F_2: {\mathcal U}_+^\omega \to  {\mathcal U}_+^\omega.$ Indeed, note that from  (A2) and (A4), we get
\begin{eqnarray*} 
0\le|F_2V(s)|&\le & \left[M\omega(s) +\beta \left(\int_S \|V\|^\theta_\omega\omega(s')^\theta q(ds'|s,a)\right)^{\frac{1}{\theta}}\right]
\\&\le&
 \left[M\omega(s)  +\beta \left( c\|V\|^\theta_\omega \omega(s)^\theta\right)^{\frac{1}{\theta}}\right] \le 
 \omega (s)\left[ M+\beta c^{\frac 1\theta}\|V\|_\omega\right].
\end{eqnarray*}
Again, by Lemma 8.5.5 in \cite{hll}, it follows that the function 
$$(s,a) \to \int_S V(s')q(ds'|s,a)$$ 
is upper semicontinuous.
Therefore, $(s,a)\to H_2(s,a,V)$ is upper semicontinuous. Consequently, by (A1) and  the Maximum Theorem of Berge 
(see Section 17.5 in \cite{ab}), the function  $F_2V$ is upper semicontinuous on $S.$ 
Now we prove that $F_2$ is $c^{\frac 1\theta}\beta$-contraction. 
Let $V_1, V_2\in {\mathcal U}^\omega_+$.  Then, by the Minkowski inequality we have
\begin{eqnarray*}
|F_2V_1(s)-F_2V_2(s)| &\le& \max_{a\in A(s)}  \beta\left| \left(\int_S V_1(s')^\theta q(ds'|s,a)\right)^{\frac 1\theta}-
\left(\int_S V_2(s')^\theta q(ds'|s,a)\right)^{\frac 1\theta}\right|\\
 &\le&
\max_{a\in A(s)}\beta
 \left(\int_S | V_1(s')-V_2(s')|^\theta q(ds'|s,a)\right)^{\frac 1\theta}\\
&\le&
\max_{a\in A(s)} \beta \left(\int_S \| V_1-V_2\|_\omega^\theta \omega(s')^\theta q(ds'|s,a)\right)^{\frac 1\theta} \\
&\le&
\beta \| V_1-V_2\|_\omega c^{\frac 1\theta} \omega(s), \quad s\in S.
\end{eqnarray*}
Now dividing all sides by $\omega(s)$ and taking the supremum over all $s\in S$  we get
$$\|F_2V_1-F_2V_2\|_\omega \le c^{\frac 1\theta} \beta  \|V_1-V_2\|_\omega.$$
Hence, $F$ is $c^{\frac 1\theta} \beta$-contraction. 
Since ${\mathcal U}_+^\omega$ is a closed subset of the Banach space ${\mathcal C}^\omega$ we conclude
from the Banach contraction principle  that there exists $w_*\in{\mathcal U}^\omega_+$ such that
$$w_*(s)=F_2w_*(s)\quad s\in S.$$
Putting $v_*:=w_*^{\frac 1{1-\rho}}$ we get that
\begin{equation}\label{be2}
v_*(s) =\max_{a\in A(s)}\left[ r(s,a)+\beta\left(\int_S  v_*(s')^{  1-\gamma}q(ds'|s,a)\right)^{\frac{1-\rho}{1-\gamma}}
\right]^{\frac{1}{ 1-\rho}},\quad s\in S.\end{equation}
for all $s\in S.$ Note that (\ref{conv2})  follows from (\ref{conv0}) with $F=F_2.$  $\Box$

\noindent{\bf Proof of Theorem \ref{thm2}}\ \ In view of Lemma \ref{l2} it remains to prove assertions $(c)$ and $(d).$
The existence of a maximiser $f_*\in \Phi$ of the right-hand side of (\ref{be2}) follows from our assumptions (A1)-(A2), (A4)
and  a measurable selection theorem, see Corollary 1 in  \cite{bp}. Again we consider the space
${\mathcal B}_+^\omega.$ 
 For any $f\in \Phi$ and $w\in{\mathcal B}_+^\omega$ we define the operator $F_{2f}$ as follows
 $$F_{2f}w(s):=H_2 (s,f(s),w),\quad s\in S.$$
Similarly as in the proof of Theorem \ref{thm1} we may show that $F_{2f}w\in {\mathcal B}_+^\omega.$ Moreover, 
$F_{2f}$ is $c^{\frac 1\theta}\beta$-contraction.
Therefore, the operator $F_{2f}$ has s a unique fixed point $w_f\in{\mathcal B}_+^\omega$ and for any function $w_0  \in{\mathcal B}_+^\omega$
$$\lim_{n\to\infty} \|F^n_{2f}w_0-w_f\|_\omega=0.$$
Observe also that $F_{2f}$ is monotone, i.e. $F_{2f}w_1\ge F_{2f}w_2$ whenever $w_1\ge w_2$ and $w_1,w_2\in{\mathcal B}_+^\omega.$  
For any $f\in F$ we have $F_{2f}w_*\le w_*.$ By iteration of this inequality we get
$F^n_{2f}w_*\le w_*$ and consequently $w_f\le w_*.$ 
On the other hand,  $F_{2f_*}w_*=w_*=F_2w_*.$ Again iterating $F_{2f_*}$ we obtain $\lim_{n\to\infty} F^n_{2f_*} w_*=w_{f_*}=w_*.$
Thus, 
$$
w_f\le w_*=w_{f_*}\quad\mbox{for all }\ f\in \Phi,
$$
 which via transformation $v:=w^{\frac 1{1-\rho}}$ with $1-\rho>0$ gives that 
$$
v_f\le v_*=v_{f_*}\quad\mbox{for all }\ f\in \Phi.
$$
This means that  $f_*$ is  an optimal stationary policy.
 $\Box$

\noindent{\bf Proof of Theorem \ref{thm3}}\ \ Proceeding as in Lemma \ref{l2} we obtain that $F_3$   is $c^{\frac 1\theta}\beta$-contraction.
Moreover,  from (B2) and (B4) for any $V\in{\mathcal C}_+^\omega$ we have
$$0<\|F_3V\|_\omega \le 
  M+\beta c^{\frac 1\theta}\|V\|_\omega.$$
Note that Lemma 8.5.5 in \cite{hll}, (B1)-(B2), (B4) and the Maximum Theorem of Berge (Section 17.5 in \cite{ab}) imply 
that $F_3V$ is continuous on $S.$
Since ${\mathcal C}_+^\omega$ is a closed subset of the Banach space ${\mathcal C}^\omega$ we conclude
from the Banach contraction principle  that there exists $w_*\in{\mathcal C}^\omega_+$ such that
$$w_*(s)=F_3w_*(s)\quad s\in S.$$
Note that $w_*(s)\ge \min_{a\in A(s)} r(s,a)>0$ for every $s\in S.$
Putting $v_*:=w_*^{\frac 1{1-\rho}}$ we get that
\begin{equation}\label{be3}
v_*(s) =\max_{a\in A(s)}\left[ r(s,a)+\beta\left(\int_S  v_*(s')^{  1-\gamma}q(ds'|s,a)\right)^{\frac{1-\rho}{1-\gamma}}
\right]^{\frac{1}{ 1-\rho}},\quad s\in S.\end{equation}
for all $s\in S.$ Clearly,  (\ref{conv3})  follows from (\ref{conv0}) with $F=F_3.$ Now, we prove (c) and (d).
Let $f_*\in \Phi$ be  a maximiser  of the right-hand side of (\ref{be3}).
 For any $f\in \Phi$ and $w\in{\mathcal B}_+^\omega$ set 
 $$F_{3f}w(s):=H_3 (s,f(s),w),\quad s\in S.$$
It follows that $F_{3f}w\in {\mathcal B}_+^\omega$ and 
$F_{3f}$ is $c^{\frac 1\theta}\beta$-contraction.
Therefore, the operator $F_{3f}$ has s a unique fixed point $w_f\in{\mathcal B}_+^\omega$ and for any function $w_0\in{\mathcal B}_+^\omega$
$$\lim_{n\to\infty} \|F^n_{3f}w_0-w_f\|_\omega=0.$$
Additionally, $F_{3f}$ is also monotone, i.e. $F_{3f}w_1\ge F_{3f}w_2$ if  $w_1\ge w_2$ and $w_1,w_2\in{\mathcal B}_+^\omega.$  
For any $f\in F$ we have $F_{3f}w_*\ge w_*.$ By iteration of this inequality we get
$F^n_{3f}w_*\ge w_*$ and consequently $w_f\ge w_*.$ 
On the other hand,  $F_{3f_*}w_*=w_*=F_3w_*.$ Again iterating $F_{3f_*}$ we obtain $\lim_{n\to\infty} F^n_{3f_*} w_*=w_{f_*}=w_*.$
Thus, 
$$
w_f\ge w_*=w_{f_*}\quad\mbox{for all }\ f\in \Phi,
$$
 which via transformation $v:=w^{\frac 1{1-\rho}}$ with $1-\rho<0$ yields that 
$$
v_f\le v_*=v_{f_*}\quad\mbox{for all }\ f\in \Phi.
$$
This means that  $f_*$ is  an optimal stationary policy.
 $\Box$


\noindent{\bf Proof of Theorem \ref{thm4}}\ \ Similarly as in Lemma \ref{l1} we show  that $F_4$ is $c\beta^{\theta}$-contraction.
Moreover,  from (B2) and (B3) for any $V\in{\mathcal C}_+^\omega$ we have
$$0< \|F_4V\|_\omega \le  \left[M+\beta c^{\frac 1\theta}\|V\|_\omega^{\frac 1\theta}\right]^\theta.
$$
Lemma 8.5.5 in \cite{hll}, (B1)-(B3)  and the Maximum Theorem of Berge (Section 17.5 in \cite{ab}) imply
the continuity of  $F_4V$ on $S.$
From the Banach contraction principle  there exists $w_*\in{\mathcal C}^\omega_+$ such that
$$w_*(s)=F_4w_*(s)\quad s\in S.$$
By (B2) $w_*(s)\ge \min_{a\in A(s)} r(s,a)>0$ for every $s\in S.$
Putting $v_*:=w_*^{\frac 1{1-\gamma}}$ we get that
\begin{equation}\label{be4}
v_*(s) =\max_{a\in A(s)}\left[ r(s,a)+\beta\left(\int_S  v_*(s')^{  1-\gamma}q(ds'|s,a)\right)^{\frac{1-\rho}{1-\gamma}}
\right]^{\frac{1}{ 1-\rho}},\quad s\in S.\end{equation}
for all $s\in S.$ Obviously,  (\ref{conv4})  follows from (\ref{conv0}) with $F=F_4.$ It remains to prove (c) and (d).
Let $f_*\in \Phi$ be  a maximiser  of the right-hand side of (\ref{be4}).
 For  $f\in \Phi$ and $w\in{\mathcal B}_+^\omega$ set 
 $$F_{4f}w(s):=H_4 (s,f(s),w),\quad s\in S.$$
 Then, $F_{4f}w\in {\mathcal B}_+^\omega$ and 
$F_{4f}$ is $c^{\frac 1\theta}\beta$-contraction.
Hence, the operator $F_{4f}$ has s a unique fixed point $w_f\in{\mathcal B}_+^\omega$ and for any function $w_0\in{\mathcal B}_+^\omega$
$$\lim_{n\to\infty} \|F^n_{4f}w_0-w_f\|_\omega=0.$$
Additionally, $F_{4f}$ is also monotone, i.e. $F_{4f}w_1\ge F_{4f}w_2$ if  $w_1\ge w_2$ and $w_1,w_2\in{\mathcal B}_+^\omega.$  
For any $f\in F$ we have $F_{4f}w_*\ge w_*.$ By iteration  we obtain
$F^n_{4f}w_*\ge w_*$ and finally $w_f\ge w_*.$ 
On the other hand,  $F_{4f_*}w_*=w_*=F_4w_*.$ Iterating $F_{4f_*}$ we get $\lim_{n\to\infty} F^n_{4f_*} w_*=w_{f_*}=w_*.$
Thus, 
$$
w_f\ge w_*=w_{f_*}\quad\mbox{for all }\ f\in \Phi,
$$
 which via transformation $v:=w^{\frac 1{1-\gamma}}$ with $1-\gamma<0$ yields that 
$$
v_f\le v_*=v_{f_*}\quad\mbox{for all }\ f\in \Phi.
$$
This proves that  $f_*$ is  an optimal stationary policy.
 $\Box$


\section{Remarks on an application of 
  Du's fixed point theorem and Banach's contraction principle}\label{sec:r}

In this section, we consider bounded utility functions. 
\subsection{Du's theorem and convex/concave operators}\label{duthm}

We first comment on some issues arising from applications
of the fixed point theorem of \cite{du}.
Recall that  $\mathcal C$ is the Banach space of all bounded continuous functions on $S$ 
endowed with the supremum norm $\|\cdot \|.$ By ${\mathcal C}_+$ we  denote the cone of all 
non-negative functions in $\mathcal C$. Let $g_1,\ g_2 \in {\mathcal C}_+$ and $g_1< g_2.$   By $I=[g_1,g_2]$ we denote the interval
in ${\mathcal C}_+$, i.e., $I=\{ v \in {\mathcal C}_+: g_1(s) \le v(s)\le g_2(s)\ \forall  s\in S\}.$ 
The following result follows from 
the fixed point theorem of   \cite{du}.

\begin{props}\label{prop1}
 Let $T:I\to I$ be a monotone increasing operator. If either (i) $T$ is convex on $I$ and for some $\epsilon \in (0,1)$
\begin{equation}
\label{t1}
Tg_2 \le  (1-\epsilon)g_2 +\epsilon g_1, 
\end{equation}
or (ii)  $T $ is concave on $I$ and  
 for some $\epsilon \in (0,1)$
\begin{equation}
\label{t2}
Tg_1\ge  (1-\epsilon)g_1 +\epsilon g_2,
\end{equation}
then $T$ has a unique fixed point $v_* \in I.$ 
Moreover, for any $v_0\in I$ and   the sequence defined by $v_n =Tv_{n-1},$ $n\in \mathbb{N},$ there exists $B>0$ such that
\begin{equation}
\label{gconv}
\|v_n - v_*\| \le B (1-\epsilon)^n,\ \ n\in \mathbb{N}. 
\end{equation} 
\end{props}

Since ${\mathcal C}_+$ is a normal cone in $\mathcal C$ with constant $N=1$, from the proofs of the theorem of   \cite{du} 
or Theorem 2.1.2 in \cite{zhang}, it follows that, if  $T$ satisfies $(i),$ then
\begin{equation}
\label{b1}
B= \|g_1-g_2\|+ \frac{2\|Tg_2 -g_2\|}{\epsilon^2}
\end{equation} 
and if $T$ satisfies $(ii),$  then
\begin{equation}
\label{b2}
B= \|g_1-g_2\|+ \frac{2\|Tg_1 -g_1\|}{\epsilon^2}.
\end{equation}

\begin{rk}\label{dubco} {\rm 
Inequalities (\ref{t1}) and (\ref{t2}) are called boundary conditions. 
They are crucial in proving the geometric convergence of iterations of the operator $T$ to a fixed point
that belongs to the interior of the set $I.$ }
\end{rk}

Below we consider two cases discussed in 
\cite{rs1} who applied the transformation $v=w^{1-\gamma}$. 
Under this transformation $H$ defined in (\ref{H}) can be written as
 $${\mathcal H}(s,a,v)=
\left[ r(s,a)+\beta\left(\int_S v(s')q(ds'|s,a)\right)^{\frac{1}{\theta}}\right]^\theta, \quad (s,a)\in\mathbb{K}.
$$  
Recall also that  $\theta =\frac{1-\gamma}{1-\rho}. $

\noindent {\bf Case 1:}  $0<\rho<\gamma<1,$ i.e., 
 $\theta\in(0,1).$ 
In this setting,
Bellman equation (\ref{be})  is equivalent to
\begin{equation}
\label{beqv}
 v(s)= Tv(s):=    \max_{a\in A(s)} {\mathcal H}(s,a,v),
\end{equation}
From the convexity of the function $\psi(t)=(b+\beta t^{1/\theta})^\theta$ on $\mathbb{R}_+$ 
for any $b\ge 0$, it follows that for each
$v_1$, $v_2$ and $\xi \in (0,1)$,
\begin{equation} \label{tconvex}
{\mathcal H}(s,a,\xi v_1+(1-\xi)v_2)\le
\xi  {\mathcal H}(s,a,v_1)+ (1-\xi ){\mathcal H}(s,a,v_2).
\end{equation}
Hence,
$$
T(\xi v_1(s) + (1- \xi) v_2)(s)  \le \xi Tv_1(s) + (1-\xi)  Tv_2(s),
$$
i.e., $T$ is convex. 

\noindent {\bf Case 3:}  $1<\rho<\gamma,$ i.e.,  $\theta >1.$  
In this setting,  
Bellman equation (\ref{be})  is equivalent to
\begin{equation}
\label{beqvmin}
 v(s)= Tv(s):=    \min_{a\in A(s)} {\mathcal H}(s,a,v), \quad s\in S.
\end{equation}
Since the function $\psi(t)=(b+\beta t^{1/\theta})^\theta$ is concave on $\mathbb{R}_+$  
for any $b\ge 0,$ it follows that for each
$v_1$, $v_2$ and $\xi \in (0,1)$,
\begin{equation} \label{tconcave}
{\mathcal H} (s,a,\xi v_1+ (1-\xi)v_2)\ge
\xi {\mathcal H}(s,a,v_1)+ (1-\xi ){\mathcal H}(s,a,v_2).
\end{equation}
Hence,
$$
T(\xi v_1(s) + (1- \xi) v_2)(s)  \ge \xi Tv_1(s) + (1-\xi)  Tv_2(s),
$$
i.e., $T$ is concave.
 
Below we comment on two other cases for which the Banach contraction
mapping theorem is applied in Subsection \ref{F}
after suitable transformations. We explain why the fixed point theorem of \cite{du}
cannot be applied.

\noindent {\bf Case 2:}  $0<\gamma< \rho<1,$ i.e., $\theta>1.$
Now Bellman equation (\ref{be})  is  equivalent to
 (\ref{beqv}). 
Now $v \to {\mathcal H}(s,a,v)$ is concave. 
Therefore, (\ref{tconcave}) holds for each
$v_1$, $v_2$ and $\xi \in (0,1).$ Hence,
$$
\max_{a\in A(s)}{\mathcal H}(s,a,\xi v_1+ (1-\xi)v_2)\ge
 \max_{a\in A(s)} ( \xi {\mathcal H}(s,a,v_1)+ (1-\xi ){\mathcal H}(s,a,v_2)),
$$
but
$$
 \max_{a\in A(s)} ( \xi {\mathcal H}(s,a,v_1)+ (1-\xi {\mathcal H}(s,a,v_2))\not\ge 
\xi \max_{a\in A(s)}  {\mathcal H} (s,a,v_1)+ (1-\xi )
\max_{a\in A(s)} {\mathcal H} (s,a,v_2).
$$
Thus, in  this case $T$  is not concave and the theorem of \cite{du} cannot be applied.  

\noindent {\bf Case 4:}   $1<\gamma< \rho,$ i.e.,  $\theta \in(0,1).$ Here,
Bellman equation (\ref{be})  is equivalent to
(\ref{beqvmin}).
The function 
$\psi(t)=(b+\beta t^{1/\theta})^\theta$ is convex on $\mathbb{R}_+$ and therefore, 
 (\ref{tconvex}) holds for each
$v_1$, $v_2$ and $\xi \in (0,1).$ 
Hence, 
$$
T(\xi v_1(s) + (1- \xi) v_2)(s)  \le 
 \min_{a\in A(s)}\left(   \xi {\mathcal H}(s,a,v_1)+ (1-\xi ){\mathcal H} (s,a,v_2)\right)
$$
but 
$$
\min_{a\in A(s)}\left(   \xi {\mathcal H} (s,a,v_1)+ (1-\xi ){\mathcal H}(s,a,v_2)\right) \not\le
\xi \min_{a\in A(s)} {\mathcal H}(s,a,v_1)+
(1-\xi )\min_{a\in A(s)} {\mathcal H} (s,a,v_2).$$
Thus, $T$ is not convex and the theorem of \cite{du}  cannot be applied to this case.

\subsection{Comparing the speed of convergence in
  Du's theorem and Banach's theorem} \label{gc}

Our aim in this section is to compare the iteration algorithms proposed by   \cite{rs1,rs2} with our modifications. 
We are interested in estimating the speed of convergence in both cases. 

We remind that from the proof of the Banach contraction mapping principle, it follows that,  
 if $F:{\mathcal C}_+ \to {\mathcal C}_+$ is a $\delta$-contraction, then $F$ has  a unique fixed point $w_*
\in {\mathcal C}_+$ and for each $n\in \mathbb{N} $
\begin{equation}
\label{b3}
\|w_n -w_*\|\le L\delta^n   \quad \mbox{where}\quad
L= \frac{\|F\hat{0} -\hat{0}\|}{1-\delta}, 
\end{equation}
$w_0\equiv \hat{0}$  and $w_n =Fw_{n-1},$
$n\in\mathbb{N}.$

\noindent
{\bf Case 1:} $0<\rho<\gamma<1.$
We assume that 
 $$ m= \min_{(s,a)\in\mathbb{K}}r(s,a)=0\quad\mbox{and}\quad 
M= \max_{(s,a)\in\mathbb{K}}r(s,a).$$ 
We fix $y>0$ 
and similarly as in \cite{rs1} we consider
the constant functions
$$g_1(s):= \hat{0}\quad\mbox{and}\quad g_2(s):
=  \left(\frac{M+y}{1-\beta}\right)^\theta,\ \  s\in S.$$
As noted by \cite{rs1}, the operator $T$ considered in this case is convex. 
Thus, the geometric convergence of iterations of $T$ with $v_0=g_1=\hat{0}$
can be analysed using the formula (\ref{b1}). Below we try to estimate the values of $B$ and $1-\epsilon$ when $\theta =1/2.$

\begin{ex}\label{ex1} {\rm  Assume that    $M=1$,  
 $\beta=0.9,$ $ \rho= 1/2 $ and $\gamma = 3/4.$ 
Then $\theta = 1/2.$ 
Let 
 $$\min_{s\in S}\max_{a\in A(s)}r(s,a)= 0. $$ 
  We have   
$$\|g_1-g_2\|= \sqrt{10(1+y)}\quad\mbox{and}\quad
Tg_2(s)= \sqrt{\max_{a\in A(s)}r(s,a)+9(1+y)}$$
and consequently,
 \begin{eqnarray*} 
\|Tg_2-g_2\| &=& \sqrt{10+10y} 
-\sqrt{9(1+y)+\min_{s\in S}\max_{a\in A(s)}r(s,a)}\\
&=&\sqrt{10+10y} - 3\sqrt{1+y}= (\sqrt{10}-3)\sqrt{1+y}.
\end{eqnarray*}
It is easy to see that   (\ref{t1})  holds if 
$$
\sqrt{1 +9(1+y)}\le (1-\epsilon)\sqrt{10+10y}.
$$
 Hence,
$$
\epsilon \le \epsilon (y):= 1-
\sqrt{\frac{ 1+0.9y}{1+y}}, \quad y>0.$$
Then 
$$1-\epsilon (y) = \sqrt{\frac{ 1+0.9y}{1+y}}.$$
To obtain the  fastest geometric convergence of iterations of the operator $T$ with
starting function $v_0=\hat{0}$, we  apply the following  method. We can  minimise  
$B(y)(1-\epsilon(y))$  where $B(y) =B$ is given in (\ref{b1}). 
Note that
$$B(y)= B =  
 \sqrt{10+10y}+  \frac{ 2(\sqrt{10}-3)\sqrt{1+y}}{\left(1-\sqrt{\frac{1+0.9y}{1+y}}\right)^2}$$ 
  and
$$
B(y)(1-\epsilon (y))=  
\left(\sqrt{10 }+  \frac{ 2(\sqrt{10}-3)}{\left(1-\sqrt{\frac{1+0.9y}{1+y}}\right)^2}\right)\sqrt{1+0.9y}
$$
where $ y>0.$ 
The minimum value of $B(y)(1-\epsilon (y))$ over the interval $(0,\infty)$ is attained at $y_0 \approx 4.10707$ and   equals
$B(y_0)(1-\epsilon (y_0)) \approx 424.197.$ 
Next $1-\epsilon (y_0)\approx 0.958947704> 
\sqrt{\beta} =\sqrt{0.9} \approx 0.948683298.$ The reason for these results is that $T$ is not  treated as a contraction mapping in 
Proposition  \ref{prop1}.
We note that $$\inf_{y>0}(1-\epsilon(y))= \lim_{y\to\infty}(1-\epsilon(y))= \sqrt{\beta}
=\sqrt{0.9},$$ but $\lim_{y\to\infty}B(y)=\infty.$ 
 Thus, minimising only the expression $1-\epsilon(y)$  makes no sense. 

The operator $F_1$ proposed in Subsection \ref{F}   
to study iterations with starting function $v_0=\hat{0}$ is as $T$ in (\ref{beqv}) with $\theta=1/2.$ 
 We have shown that $F_1$ is a $\sqrt{\beta}$-contraction. 
 Since the algorithm using the operator $F_1$ starts with $v_0=\hat{0},$ the constant 
 $L$ in (\ref{b3}) is equal to $1/(1-\sqrt{\beta})= 1/ (1-\sqrt{0.9})=19.4868.$
It is obvious that the speed of convergence of iterations of the operator $F_1$ viewed as a $\sqrt{\beta}$-contraction
is much better than in the former case.
It should be noted that the solution to Bellman equation (\ref{be}) in this example is $w_*^4.$}
\end{ex}

\noindent
{\bf  Case 3:} $1<\rho<\gamma.$
We now assume that 
  $$0<m= \min_{(s,a)\in\mathbb{K}}r(s,a)\quad\mbox{and}\quad 
M= \max_{(s,a)\in\mathbb{K}}r(s,a).$$
We fix $x\in(0,m)$ and  as in \cite{rs1} we consider
the constant functions
$$g_1(s):= \left(\frac{ x}{1-\beta}\right)^\theta\quad\mbox{and}\quad g_2(s):
=  \left(\frac{M}{1-\beta}\right)^\theta.$$
As noted by \cite{rs1}, the operator $T$ considered in this case is concave.
Therefore, the geometric convergence of iterations of $T$ with $v_0=g_1$
can be analysed using the formula  (\ref{gconv}) with $B$ defined in (\ref{b2}).
Below we try to estimate $B$ and $1-\epsilon$ in (\ref{b2}) 
when $\theta =2.$ 

\begin{ex}\label{ex2} {\rm
Assume that  $m=1,$ $M=5$ and $$\max_{s\in S}\min_{a\in A(s)}r(s,a)= 3.$$ 
Let $\beta=0.9,$ $ \rho= 1.25 $ and $\gamma = 1.5.$ Then $\theta = 2$ and  $x\in (0,1).$ Moreover, 
$$\|g_1-g_2\|= 100(25-x^2)\quad\mbox{and}\quad
Tg_1(s)= \min_{a\in A(s)} (r(s,a)+9x)^2\le (3+9x)^2 $$
for all $s\in S.$ 
Hence, $$\|Tg_1-g_1\|= (3+9x)^2-100x^2= 9+54x-19x^2.$$ 
One can easily see that (\ref{t2})  holds, if 
$$
(1+9x)^2 \ge (1-\epsilon)100x^2  +\epsilon  2500.
$$
 Thus, we have 
$$\epsilon \le \epsilon (x):= \frac{1+18x-19x^2}{100(25-x^2)}, \quad x\in (0,1),$$
and
$$1-\epsilon (x) = \frac{2499-18x-81x^2}{100(25-x^2)}.$$
To get the  best geometric convergence of iterations of the operator $T$ with
starting function $v_0=g_1$, we  can apply two approaches. First, we can try to minimise     
$B(x)(1-\epsilon(x)),$  where $B(x) =B$ is given in (\ref{b2}). We obtain
$$
B(x)= B = 100(25-x^2)   +\frac{2\cdot 100^2(25-x^2)^2(9+54x-19x^2)}{(1+18x-19x^2)^2},
$$
 and
$$
 B(x)(1-\epsilon (x))=  (2499 -18x-81x^2)R(x)$$
where
$$ R(x):= 
\left(1+ 
\frac{200(25-x^2)(9+54x-19x^2)}{(1+18x-19x^2)^2}\right),
$$ 
 for $x \in (0,1).$  
The minimum value of $B(x)(1-\epsilon (x))$ over the interval $(0,1)$ is attained at $x_0 \approx 0.378544$ and  is equal
to
$B(x_0)(1-\epsilon (x_0)) \approx 1.27143\times 10^7.$ This is a surprisingly large number.
Moreover, $1-\epsilon (x_0)\approx 0.997951789> \beta =0.9.$ The reason for such results is that $T$ is not a contraction mapping. 
Alternatively, one can find $x_1\in (0,1)$ that minimises $ 1-\epsilon (x)$    over $(0,1)$ and then calculate $B(x_1).$
Then
$$x_1 =\frac{79-8\sqrt{94}}{3} \approx 0.47904076,\quad \mbox{and}
$$
$$B(x_1)\approx 1.35198\cdot 10^7.$$
Note that $1-\epsilon (x_1)< 1-\epsilon (x_0),$ but
$B(x_1)(1-\epsilon (x_1)) > B(x_0)(1-\epsilon (x_0)).$
The operator $F_3$ proposed in Subsection \ref{F} to study value iterations is
$$F_3v(s)=\min_{a\in A(s)}\left(r(s,a) +\beta \sqrt{ \int_S v(s')^2q(ds'|s,a)}\right) $$
and is a $\beta$-contraction. Since the algorithm using the operator $F_3$ 
can start with $w_0=\hat{0},$ the constant $L$ in (\ref{b3}) is equal to $5/(1-\beta)= 50.$  
The difference between the algorithm using Du's theorem for concave and increasing operators and our proposal is significant.
One should remember that the solution to Bellman equation (\ref{be}) in this example is  the limit
$\lim_{n\to\infty} w_n^{1/(1-\rho)}= w_*^{- 4}.$ }
\end{ex}

\begin{rk} {\rm 
We would like to point out that in our method
of studying existence and uniqueness of a  solution
to the Bellman equation we do not need boundary conditions
mentioned in Remark \ref{dubco}. Therefore, we can drop the assumption made in the literature in case when 
$\gamma,\ \rho>1$ that the per-period utilities 
 $u(s,a)\ge u_0 $ for some $u_0>0$ and for all $(s,a)\in \mathbb{K}.$ However, if 
$\inf_{s\in S}\min_{a\in A(s)} u(s,a)=0,$
then $r(s,a)=(1-\beta)u(s,a)^{1-\rho}$ is unbounded and 
there should be a weight function $\omega$ satisfying
extra assumptions.
}
\end{rk}
 
\section{Conclusions} \label{sec:c}
We consider dynamic programming problems 
arising from  Markov decision processes with the Epstein-Zin preferences
involving certainty equivalents induced by power functions.
The values of parameters $\gamma$ and $\rho$ in the aggregators in Cases 1 and 2  that we deal with are
according to \cite{mm} of Thompson type. 
We have shown that in the two cases: (1) $0<\rho<\gamma<1$ and (2) $1<\rho <\gamma$ 
studied in \cite{rs1,rs2} one can use the Banach contraction mapping principle instead 
of a fixed point theorem of  \cite{du} for increasing and convex or concave operators acting on an ordered 
Banach space. The contraction mapping theorem is applied
to some   sets of non-negative functions defined on the state space, 
endowed with the standard weighted supremum metric. 
This is a much simpler way for studying the existence and uniqueness of recursive utilities for the models 
with the Epstein-Zin  preferences than the method using   
the Thompson metrics as in  \cite{mm} or the fixed point theorem of \cite{du} as in \cite{rs1,rs2}.   
Our assumptions on the per-period utility functions  are weaker than the corresponding ones in \cite{mm} and
\cite{rs1,rs2}. For example, we do not assume the boundary conditions 
assuring that the fixed point of the operator involved in the analysis
lies in the interior of its domain. 
In dynamic programming the boundary conditions are connected with inequalities
(\ref{t1}) and (\ref{t2}) in Proposition  \ref{prop1}. The constant
$\epsilon>0$ in (\ref{t1}) and (\ref{t2}) has an influence on the speed of 
convergence of the dynamic programming operator 
to a fixed point. Using the Banach theorem, we do not need to apply any 
boundary conditions and therefore, we
are able to provide  better estimation for the convergence  
of iterations of operators used in modified value iteration process.
The relevance of this result is that our method requires a smaller  number of iterations needed 
to reach a given accuracy of the solution to the Bellman equation than 
in the case of  Du's theorem.  
It seems that the third case considered by \cite{rs1} 
where $0<\rho<1<\gamma$ cannot be analysed by the Banach contraction mapping principle.
However, we wish to point out that our methodology applies 
to other two cases (1) $0<\gamma <\rho$ and (2) $1<\gamma <\rho,$ 
where the theorem of \cite{du} is inapplicable, since 
the dynamic programming operators are then neither convex nor concave.


\begin{thebibliography}{13}
\providecommand{\natexlab}[1]{#1}
\providecommand{\url}[1]{{#1}}
\providecommand{\urlprefix}{URL }
\expandafter\ifx\csname urlstyle\endcsname\relax
  \providecommand{\doi}[1]{DOI~\discretionary{}{}{}#1}\else
  \providecommand{\doi}{DOI~\discretionary{}{}{}\begingroup
  \urlstyle{rm}\Url}\fi
\providecommand{\eprint}[2][]{\url{#2}}

\bibitem[{Acemoglu (2009)}]{acemo}  Acemoglu, D., 2009. Introduction to Modern Economic Growth. Cambridge, 
Princeton Univ. Press.

 \bibitem[{ Aliprantis and Border (2006)  }]{ab}  Aliprantis, C.,    Border,  K.,    2006.   
  Infinite Dimensional Analysis: A Hitchhiker's  Guide.  New York, Springer.

\bibitem[{Asienkiewicz and  Ja{\'s}kiewicz (2017) }]{asj} 
 Asienkiewicz, H.,   Ja{\'s}kiewicz, A., 2017. A note on a new class of recursive utilities in
 Markov decision processes. Applicationes Mathematicae  (Warsaw),  44, 149-161.

\bibitem[{Backus et al. (2005) }]{backus} Backus, D., Routledge, B.,  Zin, S.,  2005. Exotic preferences for macroeconomics. 
NBER Macroeconomic Annual, 2014,  Vol. 19, pp. 319-414,  Cambridge,   MIT Press. 
 
\bibitem[{Balbus (2020)}]{balb}  Balbus, {\L.},   2020. On recursive utilities with non-affine aggregator 
and conditional certainty equivalent.  Econ.   Theory,  70,   551-577.

\bibitem[{Bansal and Yaron (2004) }]{by}   Bansal, R., Yaron, A.,    2004.  Risks for the long run: A potential resolution of
asset pricing puzzles.   J.   Finance,  59(4),  1481-1509. 

\bibitem[{Basu and   Bundick (2017) }]{basu}   Basu, S.,  Bundick, B., 2017.  
Uncertainty shocks in a model of effective demand.  Econometrica,   85,  937-958.

\bibitem[{B\"auerle and Ja\'skiewicz (2018)}]{bjjet} B\"auerle, N.,  Ja\'skiewicz, A., 2018.
Stochastic optimal growth model with risk-sensitive preferences. J.  Econom. Theory 173,  181-200. 

\bibitem[{B{\"a}uerle and Rieder (2011)}]{br} B{\"a}uerle, N., 
Rieder, U.,  2011.  Markov  Decision Processes with Applications to Finance. New York, Springer.

\bibitem[{Becker and Rinc\'on-Zapatero (2023)}]{brz} Becker, R.A., Rinc\'on-Zapatero, J.P., 2023.
Recursive utility for Thompson aggregators: least fixed Point, uniqueness, and approximation theories. 
Available at SSRN: https://ssrn.com/abstract=4456037.

\bibitem[{Bertsekas (2022)}]{bert} Bertsekas, D.P., 2022.  Abstract Dynamic Programming (3rd Ed.)
 Belmont, Massachusetts, Athena Scientific. 

\bibitem[{Bich et al. (2018)}]{bich}
 Bich, P.,  Drugeon, J.-P.,   Morhaim,  L. 2018. On temporal aggregators and dynamic programming.  Econ. Theory,  66, 
 787-817. 

\bibitem[{Blackwell (1965)}]{blackwell} Blackwell, D., 1965.
Discounted dynamic programming.   Ann. Math. Statist. 36(1),  226-235.

\bibitem[{ Bloise   and Vailakis (2018)}]{bv}   Bloise, G.,   Vailakis, Y., 2018.  Convex dynamic programming with 
(bounded) recursive utility.  J.  Econom. Theory,  173,  118-141.

\bibitem[{Bloise et al. (2024)}]{bvv}  Bloise, G.,  Le Van, C.,    Vailakis, Y.,  2024.  
Do not blame Bellman: it is Koopmans' fault.  Econometrica, 92(1),  111-140. 

\bibitem[{ Borovi{\'c}ka and   Stachurski (2020)  }]{bstach}  Borovi{\'c}ka, J.,  Stachurski, J., 2020. 
Necessary and sufficient conditions for existence  and uniqueness of recursive utilities.   J. Finance, 75(3), 1457-1493. 

\bibitem[{Brown and Purves (1973)}]{bp}  Brown, L.D.,   Purves, R., 1973.
Measurable selections of extrema.   Ann. Statist. 1,  902-912.

\bibitem[{Christensen (2022)}]{chris}   Christensen, T.M., 2022.   
Existence and uniqueness of recursive utilities without boundedness.  J. Econom. Theory,   200, 105413. 

\bibitem[{ Denardo (1967)}]{den} Denardo, E.V., 1967. 
Contraction mappings in the theory underlying dynamic programming. SIAM Rev., 9, 165-177. 
 
\bibitem[{ Du (1990)}]{du}  Du, Y., 1990.  Fixed points of increasing operators in ordered Banach spaces and
applications. Applicable Analysis,  38(2),  1-20.

\bibitem[{Dur\'an(2000)}]{duran}   Dur\'an J., 2003.  
Discounting long run average growth in stochastic dynamic programs.
Econ. Theory, 22, 395-413.  

\bibitem[{Epstein and Zin (1989)}]{ez}  Epstein, L.G., Zin S.E., 1989. 
Substitution, risk aversion, and the temporal behavior of consumption and asset returns: a theoretical framework.
Econometrica,  57,  937-969.

\bibitem[{Farhi and   Werning  (2008)}]{fw}   Farhi, E.,  Werning, T., 2008.  Optimal savings distortions 
with recursive preferences.   J. Monetary Econ.  55,    21-42.

\bibitem[{Hansen and  Sargent (1995)}]{hsar}  Hansen, L.P., Sargent, T.J., 1995. Discounted linear exponential quadratic Gaussian
control.  IEEE Trans. Autom.  Control 4,  968-971.  

\bibitem[{Hansen and   Scheinkman (2012) }]{hs} Hansen ,
L.P.,    Scheinkman,   J.A., 2012.  Recursive utility in a Markov environment with stochastic growth.    
Proc. Nat. Acad. Sci. U.S.A.,  109,  11967-11972.

\bibitem[{Hern\'andez-Lerma and Lasserre (1996)}]{hl}
 Hern\'andez-Lerma,O.,  Lasserre,  J.B., 1996.
 Discrete-Time Markov Control Processes: Basic Optimality Criteria. New York, Springer.

\bibitem[{Hern\'andez-Lerma and Lasserre (1999)}]{hll} Hern\'andez-Lerma O., Lasserre B., 1999. 
Further Topics on Discrete-Time Markov Control Processes.  New York, Springer. 

\bibitem[{ Ja{\'s}kiewicz and   Nowak (2011)}]{jndgaa}  Ja{\'s}kiewicz, A.,  Nowak, A.S., 2011.  
Stochastic games with unbounded payoffs: applications to robust control in economics. Dyn. Games  Appl., 1(2), 253-279.

\bibitem[{Kaplan and Violante (2014)}]{kv}  Kaplan, G.,   Violante,   G.L., 2014.  A model of the consumption response to
fiscal stimulus payments.  Econometrica,   82,  1199-1239. 

 \bibitem[{Kochenderfer (2015) }]{k} Kochenderfer,  M.J., 2015.   Decision Making under Uncertainty: Theory and
Application.  Cambridge,  MIT press.

\bibitem[{Koopmans (1960)}]{koop} 
 Koopmans, T.C., 1960. Stationary utility and impatience. Econometrica, 28, 287-309.

\bibitem[{Kreps and Porteus (1978)}]{kp} Kreps, D.M., Porteus, E.L., 1978. 
Temporal resolution of uncertainty and dynamic choice theory.
Econometrica,  46, 185-200. 

\bibitem[{Miao (2014)}]{miao} Miao, J., 2014.  Economic Dynamics in Discrete Time.  Cambridge, MIT Press.
 
\bibitem[{Marinacci and   Montrucchio (2010) }]{mm}  Marinacci, M.,  Montrucchio, L., 2010.  Unique solutions for stochastic recursive
utilities.  J.  Econom. Theory,   145(5),   1776-1804.
 
\bibitem[{Marinacci and   Montrucchio (2019) }]{mmtarski}  Marinacci, M.,  Montrucchio, L.,  2019. 
Unique Tarski fixed points. Math. Oper. Res., 44(4), 1174-1191.

\bibitem[{ Ozaki and  Streufert (1996)}]{os}   Ozaki,H.  Streufert, P.A., 1996. Dynamic programming for 
non-additive stochastic objectives. J. Math. Econ.,  25,  391-442. 

\bibitem[{Ren and  Stachurski  (2020) }]{rs1}  Ren, G.,   Stachurski, J., 2020. Dynamic programming with recursive
preferences: optimality and applications. (arxiv supplement) arXiv:1812.05748v4

\bibitem[{Ren and    Stachurski  (2021) }]{rs2}   Ren, G.,   Stachurski, J., 2021. 
Dynamic programming with value convexity.  Automatica 130, 109641.
 
\bibitem[{Sargent  and Stachurski (2024) }]{ss}   Sargent, T.,   Stachurski, J., 2024. Dynamic Programming, Vol. I: Foundations.  QuantEconomics.
 
\bibitem[{Sch\"al (1975)}]{schal}  Sch\"al, M., 1975. Conditions for optimality in dynamic programming and for the limit of
$n$-stage optimal policies to be optimal.  Z. Wahrsch. Verw. Gebiete, 32,  179-196.
 
\bibitem[{  Schorfheide  et al. (2018)}]{sch}  Schorfheide, F.,    Song, D.,  Yaron, A., 2018. Identifying long-run risks: A
Bayesian mixed-frequency approach.   Econometrica,   86(2),   617-654.  

\bibitem[{Skiadas (2009)}]{sk} Skiadas, C., 2009.  Asset Pricing Theory.  Princeton Series in Finance, Cambridge, Princeton  Univ. Press.

\bibitem[{Stachurski (2009)}]{stach}    Stachurski, J., 2009.  Economic Dynamics:  Theory and Computation. Cambridge,  MIT Press.

  
\bibitem[{Streufert (1990)}]{str} Streufert, P.A., 1990. Stationary recursive utility and dynamic programming 
under the assumption of biconvergence.  Rev. Econ. Studies 57(1), 79-97.
 
\bibitem[{Thompson (1963)}]{t}   Thompson, A.C., 1963. On certain contraction mappings in a partially ordered vector space.   
Proc.  Amer. Math. Soc. 14, 438-443.

\bibitem[{ Wessels (1977)}]{w}  Wessels, J., 1977.  Markov programming by successive approximations with 
respect to weighted supremum norms.  J.  Math. Anal.  Appl.,  58, 326-335.

\bibitem[{ Xing (2017)}]{x}  Xing, H., 2017.  Consumption-investment optimimization with 
Epstein-Zin utility in incomplete markets.  Finance Stoch.,  21, 227-262.

\bibitem[{Zhang (2012) }]{zhang}   Zhang, Z., 2012.   Variational, Topological, and Partial Order Methods with Their Applications. 
 New York, Springer.  
  
\end{thebibliography}
\end{document}